\documentclass[12pt,a4paper]{amsart}

\usepackage {amsfonts}
\usepackage[english]{babel}
\def\d{\delta}
\def\a{\alpha}
\def\p{\varphi}

\def\l{\lambda}
\def\s{\sigma}
\def\t{\theta}
\def\R{{\mathbb R}}
\def\N{{\mathbb N}}
\def\Z{{\mathbb Z}}

\def\bs{~\hfill\rule{7pt}{7pt}}

\DeclareMathOperator{\supp}{supp}

\newtheorem*{Th}{Theorem}
\newtheorem{Lem}{Lemma}

\begin{document}

\title{The crystalline measure that is not a Fourier Quasicrystal}

\author{S.Yu. Favorov}

\address{Serhii Favorov,
\newline\hphantom{iii}Faculty of Mathematics and Computer Science, Jagiellonian University,
\newline\hphantom{iii} Lojasiewicza 6, 30-348 Krakow, Poland
\newline\hphantom{iii} Faculty of Mathematics and Informatics, V.N.Karazin Kharkiv National University
\newline\hphantom{iii} Svobody sq., 4, Kharkiv, Ukraine 61022}
\email{sfavorov@gmail.com}

\maketitle {\small
\begin{quote}
\noindent{\bf Abstract.}
We construct a crystalline measure on the real line, which is not a Fourier Quasicrystal.
\medskip


\end{quote}
}

\bigskip
A complex measure $\mu$ with discrete support is called crystalline if $\mu$ is a temperate distribution, and its Fourier transform
in sense of distributions $\hat\mu$ is also a measure with discrete support.  This measure is a Fourier Quasicrystal if  the measures $|\mu|$ and $|\hat\mu|$ are also temperate distributions.
These definitions are valid both for measures on real axis and for measures on arbitrary Euclidean space.

Here we say that a set $E$ is discrete if an intersection of $E$ with any ball is a finite set.  Also, $|\nu|(E)$ means the variation of the complex measure $\nu$ on the set $E$, and $|\nu|$ means the
corresponding positive measure.

 Fourier Quasicrystals  are used as mathematical models of certain atomic structures, having a discrete diffraction pattern.
A number of papers has appeared, in which the properties of Fourier Quasicrystals and crystalline measures are studied.  Conditions for crystalline measures to be periodic are found, and nontrivial examples of Fourier Quasicrystals  are constructed (\cite{F}--\cite{Q}).

In our article we present the following result:
\begin{Th}
There is a crystalline measure on $\R$ that is not a Fourier Quasicrystal.
\end{Th}
Note that if a crystalline measure $\mu$ has the properties
$$
 |x-x'|\ge c\min\{|x|^{-h},1\}\quad\forall\,x,\,x'\in\supp\mu,\quad |y-y'|\ge c\min\{|y|^{-h},1\}\quad\forall\,y,\,y'\in\supp\hat\mu,
$$
with some $c,\,h>0$, then $\mu$ is a Fourier Quasicrystal (\cite{F0}). Therefore  points of support of the measure constructed in the theorem should  approach very quickly near infinity.

\medskip
Before constructing this example, recall the definition of the Fourier transform in the sense of distributions.

 By $S(\R)$ denote the Schwartz space of test functions $\p\in C^\infty(\R)$ with the finite norms
 $$
  N_{n,m}(\p)=\sup_{\R}\max_{k\le m} |(1+|x|^n)\p^{(k)}(x)|,\quad n,m=0,1,2,\dots
$$
 These norms generate the topology of the projective limit on $S(\R)$.  Elements of the space $S^*(\R)$ of continuous linear functionals on $S(\R)$ are called temperate distributions.
For each temperate distribution $f$ there are $C<\infty$ and $n,\,m\in\N\cup\{0\}$ such that for all $\p\in S(\R^d)$
\begin{equation}\label{N}
                           |f(\p)|\le CN_{n,m}(\p).
\end{equation}
Conversely, if a linear functional $f$ on $S(\R)$ satisfies this estimate for at least one pair $n,\,m$, then $f\in S^*(\R^d)$.

The Fourier transform of a temperate distribution $f$ is given by the equality
$$
\hat f(\p)=f(\hat\p)\quad\mbox{for all}\quad\p\in S(\R),
$$
where
$$
   \hat\p(t)=\int_{\R^d}\p(x)e^{-2\pi i xt}dx
 $$
is the Fourier transform of the function $\p$. By $\check\p$ we  denote the inverse Fourier transform of $\p$.
The Fourier transform is a bijection of $S(\R)$ onto itself and a bijection of $S^*(\R)$ onto itself.

Here and below $\d_x$ is the unit mass at the point $x$, and $\# E$ is the number of points of the finite set $E$.

\medskip

The proof of the theorem is based on the following lemmas:

\begin{Lem}[cf.\cite{F}]\label{L1}
Let $\mu$ be a discrete measure and $|\mu|\in S^*(\R)$. Then for any  $\p\in S(\R^d)$ the function $\hat\mu\star\p(t)$ is bounded.
\end{Lem}
{\bf Proof}. It is easy to prove that for any positive measure $\nu\in S^*(\R)$ there is $T<\infty$ such that $\nu(-r,r)=O(r^T)$ as $r\to\infty$ (cf.\cite{F1}).
Therefore if $\mu=\sum_\l c_\l\d_\l$, then
$$
M(r):=\sum_{-r<\l<r}|c_\l|=|\mu|(-r,r)=O(r^T)\quad\text{as}\quad r\to\infty.
$$
For any $\p\in S(\R)$ we have
$$
  \hat\mu\star\p(t)=(\hat\mu(y),\p(t-y))=(\mu(x),\check\p(x)e^{-2\pi itx})=\sum_{\l\in\supp\mu} c_\l\check\p(\l)e^{-2\pi it\l}.
$$
 Since $|\check\p(x)|\le N_{T+1,0}(\check\p)|x|^{-T-1}$ for $|x|>1$, we obtain
$$
\phantom{XXXXXXXX}  \sum_{\l\in\supp\mu}|c_\l||\check\p(\l)|\le C_0+C_1\int_1^\infty r^{-T-1}M(dr)<\infty. \phantom{XXXXXXXX}\bs
$$

\begin{Lem}\label{L2}
Let $\l_n,\,\tau_n$ be two sequences of positive numbers such that $\tau_n\to0$, $\l_n\to\infty$, and
\begin{equation}\label{g}
 \log\tau_n/\log\l_n\to-\infty\quad\mbox{as}\quad n\to\infty.
\end{equation}
Let $\mu\in S^*(\R)$ be a measure  such that its restriction to $(\l_n-1/(2\l_n),\l_n+1/(2\l_n))$ equals $\tau_n^{-2/3}(\d_{\l_n+\tau_n}-\d_{\l_n})$.
Then there is $\p\in S(\R)$  such that $\hat\mu\star\p(t)$ is unbounded.
\end{Lem}

{\bf Proof}.
By thinning out the sequence $\tau_n$, we can assume that for all $n$
\begin{equation}\label{s1}
\sum_{p<n}\tau_p^{-1/3}<(1/3)\tau_n^{-1/3},
\end{equation}
and
\begin{equation}\label{s2}
\sum_{p>n}\tau_p^{2/3}<2\tau_n ^{2/3}/(3\pi).
\end{equation}
 Let $\eta$ be a non-negative even function such that
 $$
 \eta(x)\in C^\infty(\R),\quad\eta(x)=0\mbox{ for }|x|>1/2,\quad \eta(x)=1\mbox{ for }|x|\le1/3.
 $$
Set
 $$
 \psi(x)=\sum_n \tau_n^{1/3}\eta(\l_n(x-\l_n)).
 $$
 By \eqref{g}, $\tau_n=o(1/|\l_n|^T)$ as $n\to\infty$ for every $T<\infty$. Therefore, for all $k>0,\ N\in\N$
we have
$$
\psi^{(k)}(x)=o(\l_n^{-N})\quad\text{for}\quad |x-\l_n|<1/(2\l_n).
$$
 Hence, $\psi^{(k)}(x)(1+|x|^N)$ is bounded on $\R$ for all $N$ and $k$, i.e., $\psi\in S(\R)$.
  Since  $\psi(x)=0$ for $|x-\l_n|>1/(2\l_n))$, we have for every $t\in\R$
 $$
(\hat\mu(y),\hat\psi(t-y))=(\mu(x),\psi(x)e^{-2\pi ixt})=\sum_{n=1}^\infty\tau_n^{-1/3}[\eta(\tau_n\l_n)e^{-2\pi i(\l_n+\tau_n)t}-
 \eta(0)e^{-2\pi i \l_nt}].
 $$
 For $n\ge n_0$ we have $\tau_n<1/(3\l_n)$,  therefore, $\eta(\tau_n\l_n)=\eta(0)=1$. Besides, for $t=1/(2\tau_n)$
 $$
 |e^{-2\pi i(\l_n+\tau_n)t}-e^{-2\pi i \l_nt}|=|e^{-2\pi i\tau_nt}-1|=2.
  $$
  Therefore,
 \begin{equation}\label{s3}
  |\hat\mu\star\hat\psi(t)|\ge2\tau_n^{-1/3}-2\sum_{p<n}\tau_p^{-1/3}-\sum_{p>n}\tau_p^{-1/3}|e^{-2\pi i\tau_pt}-1|\quad\text{for}\quad n\ge n_0.
 \end{equation}
Taking into account \eqref{s1}, \eqref{s2}, and the estimates
$$
|e^{-2\pi i\tau_pt}-1|\le 2\pi\tau_p t=\pi\tau_p\tau_n^{-1},
$$
we obtain that \eqref{s3} is more than $2\tau_n^{-1/3}/3$, and the convolution $\hat\mu\star\p(t)$ with $\p=\hat\psi$ is unbounded. \bs
\medskip

\begin{Lem}[Y.Meyer \cite{M}, Lemma 7, also M.N.Kolountzakis, \cite{K}]\label{L3}
Let $\a\in(0,1/6)$. For every integer $M>M_\a$ there exists an
$M$-periodic discrete measure $\s=\s_M$ such that
$$
\supp\s_M\cup\supp\hat\s_M\subset\L_M:=M^{-1}\Z\setminus[-\a M,\a M].
$$
\end{Lem}

\begin{Lem}\label{L4}
Let the measure from Lemma \ref{L3} has the form
\begin{equation}\label{a}
\s_M=\sum_{k\in\Z}\sum_{j=0}^{M^2-1} c_j\d_{kM+j/M}.
\end{equation}
 Denote by $\s^h$ the shift of $\s$ along  $h$. Then  for any $\phi\in S(\R)$ and $\tau\in(0,1),\,|h|<M/3,\,M>2$
$$
  |(\s_M^{h+\tau}-\s_M^h,\p)|\le CM^2N_{2,1}(\phi)\max_j|c_j|\tau,
$$
where $N_{2,1}(\phi)$ is defined in \eqref{N}, and $C$ is an absolute constant.
\end{Lem}
{\bf Proof}. Clearly, for all $j,\,k$ there is $\t\in(0,1)$ such that
$$
(\d_{kM+j/M+h+\tau}-\d_{kM+j/M+h},\phi)=\tau\phi'(kM+j/M+h+\t).
$$
Using the definition of $N_{2,1}(\phi)$, we obtain
$$
   |(\s_M^{h+\tau}-\s_M^h,\p)|\le \sum_{|k|\le2}\sum_{j=0}^{M^2-1}\frac{|c_j|N_{2,1}(\phi)\tau}{\{1+|kM+j/M+h+\t|^2\}}
  $$
$$
+\sum_{|k|>2}\frac{M^2\max_j|c_j|N_{2,1}(\phi)\tau}{(|k|M-M-|h|-\t|)^2}\le(5M^2+\pi^2/3)\max_j|c_j|N_{2,1}(\phi)\tau.
$$
\bs

{\bf Proof of the theorem}. Set $M_n=32^n$. Let $\tau_n$ be any sequence such that $0<\tau_n<1/(4M_n)$ and
\begin{equation}\label{t}
\log\tau_n/n\to-\infty\quad\text{as}\quad n\to\infty.
\end{equation}
We prove that with suitable choice of $h_n\in (-M_n/32,M_n/32)$ the measure
$$
  \mu=\sum_{n>n_0}\tau_n^{-2/3}(\s_{M_n}^{h_n+\tau_n}-\s_{M_n}^{h_n})
$$
satisfies the statement of the theorem. Here $\s_{M_n}$ is the measure from Lemma \ref{L3}  such that $\max_j|c_j|=1$ in representation \eqref{a}. We also suppose that  we have $c_{j'}=1$ for some  $j'$.
The number $n_0>2$ is such that  the statement of Lemma \ref{L3} is satisfied by  $\a=1/8$ and $M=M_n$ for $n>n_0$.

Applying Lemma \ref{L4} to the measure $\s_{M_n}^{h_n}$, we get  for  $\phi\in S(\R)$
$$
  |(\mu,\phi)|\le \sum_{n>n_0} C\tau_n^{1/3}M_n^2N_{2,1}(\phi).
$$
By \eqref{t}, the series $\tau_n^{1/3}M_n^2$ converges. Therefore, $\mu$ satisfies \eqref{N} and $\mu\in S(\R)$.

Since $|h_n+\tau_n|<M_n/16$, we see that for every $n$ the support of the measure $\s_{M_n}^{h_n+\tau_n}-\s_{M_n}^{h_n}$ does not intersects with $[-M_n/16,M_n/16]$,
so every bounded interval contains a finite number of points  $\supp\mu$, and the same is true for $\supp\hat\mu$. Hence, $\mu$ is a crystalline measure.

Furthermore, show that we can take $h_n$ such that $|h_n|<M_n/32$ and then $\l_n\ge M_n$ such that
\begin{equation}\label{i}
\left.\s_{M_n}^{h_n+\tau_n}-\s_{M_n}^{h_n}\right|_{(\l_n-1/(2\l_n),\l_n+1/(2\l_n))}=\d_{\l_n+\tau_n}-\d_{\l_n},
\end{equation}
 and  for all $p\neq n$
\begin{equation}\label{i1}
  \supp (\s_{M_p}^{h_p+\tau_p}-\s_{M_p}^{h_p})\cap(\l_n-1/(2\l_n),\l_n+1/(2\l_n))=\emptyset.
\end{equation}
 Set
$$
   I_n=[M_n,2M_n)\cap(M_n+j'/M_n-M_n/32,M_n+j'/M_n+M_n/32),
$$
where $j'$ is the number defined above, and
$$
I_{n,j}=[M_n+j/M_n,M_n+(j+1)/M_n),\quad 0\le j<M_n^2.
$$
Since $2M_n\le M_p/16$ for $p>n$, we get
\begin{equation}\label{i2}
I_{n,j}\cap[\supp\s_{M_p}^{h_p+\tau_p}\cup\supp\s_{M_p}^{h_p}]=\emptyset \quad\forall j,\,\forall p>n.
\end{equation}
Next, let $h_p,\,p<n$,  have already been chosen. Since for $p<n$
$$
\#\{kM_p+q/M_p\in I_n:\,k\in\N,\,q\in\N\cup\{0\}\} =M_pM_n/8,
$$
 we get
$$
 \#\{\supp\s_{M_p}^{h_p}\cap I_n\}\le M_pM_n/8,\quad  \#\{\supp\s_{M_p}^{h_p+\tau_p}\cap I_n\}\le M_pM_n/8.
$$
 Summing over $n_0<p<n$, we get
 $$
  \#\{\cup_{n_0<p<n}[\supp\s_{M_p}^{h_p+\tau_p}\cup\supp\s_{M_p}^{h_p}]\cap I_n\}<M_n^2/124.
$$
On the other hand,
$$
\#\{j:\,I_{n,j}\subset I_n\}\ge M_n^2/32-1>M_n^2/124.
$$
Therefore there exists a half-interval $I_{n,j''}\subset I_n$ such that
\begin{equation}\label{i3}
I_{n,j''}\cap(\supp\s_{M_p}^{h_p+\tau_p}\cup\supp\s_{M_p}^{h_p})=\emptyset\quad\forall\, p<n.
\end{equation}
Since $I_{n,j''}\subset I_n$, we get $|j'/M_n-j''/M_n|<M_n/32$. Taking into account that $j'-j''\in\Z$, we also get
$$
 |j'/M_n-j''/M_n|\le M_n/32-1/M_n,
$$
and for
$$
h_n:=j''/M_n-j'/M_n+1/(2M_n)
$$
we obtain $|h_n|<M_n/32$.
Set
$$
\l_n=M_n+j''/M_n+1/(2M_n).
$$
Combining \eqref{i2}, \eqref{i3}, and the embedding $(\l_n-1/(2\l_n),\l_n+1/(2\l_n))\subset I_{n,j''}$, we get \eqref{i1}.
Taking into account the estimate $\tau_n<1/(4M_n)<1/(2\l_n)$, we obtain \eqref{i}.
By \eqref{t}, we get condition \eqref{g}. Applying Lemma \ref{L2}, we find a function $\p\in S(\R)$ for which  the convolution
$\hat\mu\star\p(t)$ is unbounded on $\R$. Hence Lemma \ref{L1} implies that $|\mu|\not\in S^*(\R)$, and $\mu$ is not a Fourier Quasicrystal.   \bs
\medskip

{\bf Remark}. Y.Meyer formulated  a theorem in \cite{M} that any crystalline measure is an almost periodic distribution, i.e., its convolution with $C^\infty$-function with compact support
is an almost periodic function in the sense of Bohr. Then he wrote in \cite{M1} that his proof of this theorem is incorrect and formulated the corresponding result as Conjecture 2.1.
We have just proved that $\hat\mu$ is a crystalline measure and the convolution $\hat\mu\star\p(t)$ is unbounded on $\R$, hence this convolution is not almost periodic.
But this is only a partial refutation of the Conjecture, since the support of the function $\p$ is not a compact set.
\bigskip

 I thank the Department of Mathematics and Computer Science
of the Jagiellonian University for its hospitality and Professor Lukasz Kosinski for his interest in my work and useful discussions.

\end{document}